# Event-triggered and self-triggered stabilization of distributed networked control systems

Romain Postoyan, Paulo Tabuada, Dragan Nešić, Adolfo Anta

*Abstract*— Event-triggered and self-triggered control have recently been proposed as implementation strategies that considerably reduce the resources required for control. Although most of the work so far has focused on closing a single control loop, some researchers have started to investigate how these new implementation strategies can be applied when closing multiple-feedback loops in the presence of physically distributed sensors and actuators. In this paper, we consider a scenario where the distributed sensors, actuators, and controllers communicate via a shared wired channel. We use our recent prescriptive framework for the event-triggered control of nonlinear systems to develop novel policies suitable for the considered distributed scenario. Afterwards, we explain how self-triggering rules can be deduced from the developed event-triggered strategies.

## I. Introduction

Today's control systems are frequently implemented over networks as these types of structures present many advantages in terms of flexibility and cost. In this setup, controllers communicate with sensors and actuators through the network, not in a continuous fashion but rather at discrete time instants when the channel is available for the control system. Traditionally, the time interval between two successive transmissions is constrained to be less than a fixed constant $T$, which is called the *maximum allowable transmission interval* (MATI) (see *e.g.* [6], [14], [21]). In order to achieve a desired performance, $T$ is generally chosen as *small* as technology and network load permit. This strategy, although easy to implement and analyze, represents a conservative solution that may unnecessarily overload the communication channel. Indeed, one would expect that the transmission instants should not satisfy a prefixed bound but rather be based on the current state of the system, the channel occupancy and the desired performance. Drawing intuition from this idea, event-triggered control has been developed to reduce the need for feedback while guaranteeing satisfactory levels of performance. It involves closing the loop whenever a predefined state-dependent triggering condition is satisfied, *e.g.* [3], [4], [9], [20]. This technique reduces resource usage such as communication bandwidth or computational time and provides a high degree of robustness since the state is continuously monitored. Most work in this direction has focused so far on the so-called one packet transmission problem, which corresponds to the case when all the states are sent together in a single packet. This generally implies the collocation of all sensors and, for multiple-input control systems, collocation of all actuators as well. Such an assumption does not hold in many cases. Thus, significant work on distributed event-triggered control has recently appeared in [13], [23]. These papers focus on a setup where the sensors decide locally when they need to transmit their measurements. Proposed solutions may however be conservative as they allow several sensors to transmit simultaneously.

In this paper, we consider a network setup where a central coordinator is available and only one group of sensors or actuators (that forms a *node*) can transmit at each transmission instant. This central coordinator grants access to the node selected according to a given scheduling protocol whenever a predefined triggering rule is satisfied. We follow the same approach as in [17], which is a particular case of this study since we considered networks that have one node only. Modeling the problem using the hybrid formalism of [8], we apply the prescriptive framework in [17] to synthesize event-triggering rules for *networked control systems* (NCS). It has to be noted that existing strategies are not directly applicable since only a subset of sensors and actuators get access to the network at a transmission instant. Hence, we adapt a policy developed in [17] and develop new event-triggering rules for classes of NCS governed by uniformly globally asymptotically stable (UGAS) protocols (see [15]). Examples of UGAS protocols are the round-robin protocol or the try-once-discard protocol [21] as shown in [14]. To the best of the authors' knowledge, we believe that this is the first time that event-triggered control is addressed for such NCS.

The scheme we introduce in this paper can be seen as a centralized event-triggering policy for physically distributed control systems, which in practice would require constant communication between the sensors/actuators and the coordinator (since the decision depends on the current value of the system states). To overcome this problem, we show how self-triggered implementations (see [1], [2], [22]) can be derived from a known event-triggering strategy by applying the techniques in [2]. Under this paradigm, the coordinator decides the next time instant at which communication should occur based on the last received information.

The remainder of the paper is organized as follows. The problem we solve in this paper is described in Section II. In Section III we recall the hybrid Lyapunov-based framework proposed in [17] to analyse the stability of NCS. This framework is used as a baseline to derive event-triggering conditions in Section IV. Section V describes how such

This technical report is an extended version of [18]. This work was supported by the Australian Research Council under the Future Fellowship and Discovery Grants Schemes, the NSF awards 0820061 and 0834771 and the Alexander von Humboldt Foundation.

R. Postoyan is with the Centre de Recherche en Automatique de Nancy, UMR 7039, Nancy-Université, CNRS, France `romain.postoyan@cran.uhp-nancy.fr`

P. Tabuada is with the Department of Electrical Engineering, University of California at Los Angeles, Los Angeles, CA 90095-1594, USA `tabuada@ee.ucla.edu`

D. Nešić is with the Department of Electrical and Electronic Engineering, the University of Melbourne, Parkville, VIC 3010, Australia `d.nesic@ee.unimelb.edu.au`

A. Anta is with the Technische Universität Berlin & Max Planck Institute für Dynamik komplexer technischer Systeme, Germany `anta@control.tu-berlin.de`

event-triggered conditions can be emulated by means of the self-triggered formulation. Finally, an illustrative example of a distributed system is covered in Section VI to show the benefits of this strategy. The proofs are given in the Appendix.

## II. PROBLEM STATEMENT

We omit the description of the notation and definitions used throughout this paper and refer the reader to Section II in [17]. Consider the plant:

$$\dot{x}_P = f_P(x_P, u), \quad (1)$$

where $x_P \in \mathbb{R}^{n_P}$ is the plant state and $u \in \mathbb{R}^{n_u}$ the control input. The following stabilizing dynamic state-feedback controller is designed:

$$\dot{x}_C = f_C(x_C, x_P), \quad u = g_C(x_C, x_P), \quad (2)$$

where $x_C \in \mathbb{R}^{n_C}$ is the controller state. We consider a scenario where the controller (2) communicates with the plant (1) via a shared centralized network.

Since state measurements and control inputs are no longer continuously available but only at given transmission instants $t_j$, $j \in \mathbb{Z}_{>0}$, systems (1) and (2) become:

$$\begin{array}{rcll} \dot{x}_P & = & f_P(x_P, \hat{u}) & \forall t \in [t_{j-1}, t_j] \\ \dot{x}_C & = & f_C(x_C, \hat{x}_P) & \forall t \in [t_{j-1}, t_j] \\ u & = & g_C(x_C, \hat{x}_P), & \end{array} \quad (3)$$

where $\hat{x}_P$ and $\hat{u}$ denote the variables respectively generated from the most recent transmitted plant state and control input through the network. Between two transmission instants, they are generated by the in-network processing algorithm modeled by functions $\hat{f}_P$ and $\hat{f}_C$:

$$\begin{array}{rcll} \dot{\hat{x}}_P & = & \hat{f}_P(x_P, x_C, \hat{x}_P, \hat{u}) & \forall t \in [t_{j-1}, t_j] \\ \dot{\hat{u}} & = & \hat{f}_C(x_P, x_C, \hat{x}_P, \hat{u}) & \forall t \in [t_{j-1}, t_j]. \end{array} \quad (4)$$

Zero-order-hold devices are often used so that $\hat{x}_P$ and $\hat{u}$ are kept constant on $[t_{j-1}, t_j]$ i.e. $\hat{f}_P = 0$ and $\hat{f}_C = 0$. Nevertheless, we allow for other types of implementations.

Sensors and actuators are grouped into $l$ nodes depending on their spatial location. At each transmission instant, a single node gets access to the channel according to the scheduling protocol and transmit its data. We model this process as follows:

$$\begin{array}{rcl} \hat{x}_P(t_j^+) & = & x_P(t_j) + h_P(j, e(t_j)) \\ \hat{u}(t_j^+) & = & u(t_j) + h_u(j, e(t_j)), \end{array} \quad (5)$$

where $e = (e_{x_P}, e_u)$, with $e_{x_P} = \hat{x}_P - x_P$ and $e_u = \hat{u} - u$, denoting the networked-induced error. Thus, vector $e$ is partitioned as $e = (e_1, \ldots, e_l)$. At each transmission instant $t_j$, functions $h_P, h_u$ are typically such that if the node $i$ gets access to the network, the corresponding error $e_i$ experiences a jump while the other components of $e$ remain unchanged; usually $e_i(t_j^+) = 0$ but this assumption is not needed in general. In that way, $h_P$ and $h_C$ can be used to model common scheduling protocols such as round-robin (RR) or try-once-discard (TOD) [21], see [14] for more details.

The sequence of transmission instants $t_j$, $j \in \mathbb{Z}_{>0}$, is traditionally defined such that $\varepsilon \leq t_j - t_{j-1} \leq T$ (as in [21], [14], [6]), where $T \in \mathbb{R}_{>0}$ denotes the MATI and $\varepsilon \in \mathbb{R}_{>0}$ is an arbitrary small constant which models the fact that there exists a minimum amount of time between two transmissions. In this study, we resort to a different paradigm: transmissions are triggered according to a criterion that depends on the variables of the overall system. The underlying idea behind this paradigm is to reduce the usage of the communication bandwidth by transmitting data only when needed to ensure the desired stability properties.

We consider the following implementation architecture. Sensors, actuators and controllers exchange information through a network, where the schedule is dynamically decided by a central coordinator. The order at which the network is assigned to each node is defined by means of the protocol. The time instants at which communication needs to be established are decided by the central coordinator according to the current state of the dynamical system, in the spirit of event-triggered control. This central node receives information from the sensors and the controller and evaluates a so-called triggering rule in order to decide whether communication is needed to guarantee stability for the control system. Since the triggering condition usually depends on the state of the plant and the controller, such setup requires in general continuous communication between the sensors, the controller and the central coordinator. To overcome this handicap, we propose in this paper a self-triggered implementation that emulates the designed event-triggered policy, where the next transmission is decided based on the last data received by the coordinator (see Section V). Since transmission times are known in advance under this policy, the self-triggered policy also facilitates the schedulability analysis for the network.

We model the problem using the hybrid formalism of [8], similar to [6], [7]. We group together the states of the plant and the controller in the variable $x = (x_P, x_C) \in \mathbb{R}^{n_x}$ and we denote by $\kappa \in \mathbb{Z}_{\geq 0}$ the counter variable that may be required to model protocols such as round-robin for instance (see [14]). It has to be noted that additional variables may also be introduced for designing the triggering rule. For instance, we will see in Section IV-A that the event-triggering strategy in [20] is not applicable to the considered distributed NCS unless we introduce an appropriate auxiliary variable. We will also show in Section IV-B that the time-triggered policy in [6] can be modified to obtain event-triggering rules thanks to the use of a clock-like variable. Thus, we denote by $\eta \in \mathbb{R}^{n_\eta}$ all auxiliary variables. The model can be written as:

$$\left.\begin{array}{rcl} \dot{x} & = & f_x(x, e) \\ \dot{e} & = & f_e(x, e) \\ \dot{\kappa} & = & 0 \\ \dot{\eta} & = & f_\eta(x, e, \kappa, \eta) \end{array}\right\} \quad q \in C \\ \left.\begin{array}{rcl} x^+ & = & x \\ e^+ & = & h_e(\kappa, e) \\ \kappa^+ & = & \kappa + 1 \\ \eta^+ & = & h_\eta(x, e, \kappa, \eta) \end{array}\right\} \quad q \in D, \quad (6)$$

where $q = (x, e, \kappa, \eta)$. We use $\dot{q} = f_q(q)$ and $q^+ = h_q(q)$ to denote (6). The sets $C$ and $D$ are closed, included in $\mathbb{R}^{n_q}$ ($n_q = n_x + n_e + 1 + n_\eta$) and respectively denote the flow and the jump set. Typically, the system flows on $C$ and experiences a jump on $D$, where the triggering condition is satisfied. When $q \in C \cup D$, the system can either jump or flow, the latter only if flowing keeps $q$ in $C$. We call $e^+ = h_e(\kappa, e)$ the protocol where $h_e = (h_P, h_u)$ as in [14],

[6]. Functions $f_\eta$, $h_\eta$,

$$f_x : (x, e) \mapsto \begin{pmatrix} f_P(x_P, g_C(x_C, \hat{x}_P) + e_u) \\ f_C(x_C, \hat{x}_P) \end{pmatrix}$$

$$f_e : (x, e) \mapsto \begin{pmatrix} \hat{f}_P(x_P, x_C, \hat{x}_P, g_C(x_C, \hat{x}_P) + e_u) \\ -f_P(x_P, g_C(x_C, \hat{x}_P) + e_u) \\ \hat{f}_C(x_P, x_C, \hat{x}_P, g_C(x_C, \hat{x}_P) + e_u) \\ -\frac{\partial g_C}{\partial x_C}(x_C, \hat{x}_P) f_C(x_C, \hat{x}_P) - \frac{\partial g_C}{\partial \hat{x}_P}(x_C, \hat{x}_P) \\ \times f_P(x_P, x_C, \hat{x}_P, g_C(x_C, \hat{x}_P) + e_u) \end{pmatrix} \quad (7)$$

where $\hat{x}_P = x_P + e_{x_P}$, are assumed to be continuous.

The main problem addressed in this paper is to define appropriate event-triggering rules, that is, to define appropriate flow and jump sets $C$ and $D$ for system (6) in order to ensure asymptotic stability properties of (6) while reducing the number of transmissions as much as possible. Afterwards, we explain how self-triggering conditions may be derived from a known event-triggering criterion.

### III. A PRESCRIPTIVE FRAMEWORK FOR THE EVENT-TRIGGERED CONTROL OF NCS

We recall in this section the framework of [17], originally developed for sampled-data systems. It is based on the following theorem that provides sufficient conditions that ensure asymptotic stability properties for system (6). It can be regarded as a variation of the general results in [5].

**Theorem 1.** *Consider system (6) and suppose $h_q(D) \subset (C \cup D)$ and that there exist a locally Lipschitz function $R : \mathbb{R}^{n_q} \to \mathbb{R}$ and a continuous function $v : \mathbb{R}^{n_\eta} \to \mathbb{R}^{n_v}$ with $n_v \leq n_\eta$ such that the following conditions hold:*

(i) *There exist $\underline{\alpha}_R, \overline{\alpha}_R \in \mathcal{K}_\infty$ such that for any $q \in C \cup D$: $\underline{\alpha}_R(|(x, e, v(\eta))|) \leq R(q) \leq \overline{\alpha}_R(|(x, e, v(\eta))|).$*
(ii) *There exists $\alpha_R \in \mathcal{K}_\infty$ such that for all[1] $q \in C$: $R^\circ(q; f_q(q)) \leq -\alpha_R(R(q)).$*
(iii) *For all $q \in D$, $R(h_q(q)) \leq R(q)$.*
(iv) *Solutions to (6) have a semiglobal dwell time[2] on $\mathbb{R}^{n_q} \backslash \mathcal{A}$, where $\mathcal{A} = \{q : (x, e, v(\eta)) = 0\}$.*

*Then the set $\mathcal{A}$ is S-GAS.*

Theorem 1 can be used as a framework for the synthesis of event-triggering rules for (6). The main idea is to design the triggering criterion so that there exists a Lyapunov function for the overall system (6) that decreases on flows, does not increase at jumps and guarantees the existence of a minimal interval of times between two jumps outside the stable set. This approach has been used to investigate the stability of other types of hybrid systems (*e.g.*, see [16], [6]). General guidelines on how to apply Theorem 1 to synthesize triggering conditions for system (6) can be found in Section IV in [17]. The main difference with [17] is that the $e$-dependency of the Lyapunov function $R$ will depend on the considered scheduling protocol as we show it in Section IV. The following result is used to verify the existence of semiglobal dwell times in Section IV. It is a corollary of Lemma 1 in [17].

**Lemma 1.** *Consider system (6) and assume that $h_q(D) \subset (C \cup D)$ and items (i)-(iii) of Theorem 1 are satisfied. If*

---
[1] We consider the $R^\circ(q; f_q(q))$, the Clarke derivative of $R$ (see [17]), by abuse of notation, although $R$ is not necessarily locally Lipschitz in $\kappa$. This is justified since the component of $f_q(q)$ corresponding to $\kappa$ is 0.
[2] See Definition 2 in [17]

*for any $\mu \in \mathbb{R}_{>0}$ there exists a locally Lipschitz function $\psi : \Theta(\mu) \to \mathbb{R}_{\geq 0}$ in $(x, e, \eta)$, where $\Theta(\mu) = \{q \in C \cup D : R(q) \in (0, \mu]\}$ such that:*

(i) *There exists $a \in \mathbb{R}_{\geq 0}$ such that for any $q \in D$ with $h_q(q) \in \Theta(\mu)$: $\psi(h_q(q)) \leq a$.*
(ii) *There exists $b > a$ such that for any $t_j \leq t$ with $(t, j) \in \text{dom } \phi$: $(\psi(\phi(t, j)) < b) \Rightarrow (\phi(t, j) \in C \backslash D)$.*
(iii) *There exists a continuous non-decreasing function $\lambda : \mathbb{R}_{>0} \to \mathbb{R}_{\geq 0}$ such that for all $q \in \Theta(\mu) : \psi^\circ(q; f_q(q)) \leq \lambda(\psi(q)).$*

*Then solutions to (6) have a semiglobal dwell time on $\mathbb{R}^{n_q} \backslash \mathcal{A}$ where $\mathcal{A} = \{q : (x, e, v(\eta)) = 0\}$.*

### IV. EVENT-TRIGGERED STRATEGIES

We apply the framework of Section III to synthesize event-triggering rules for NCS. Two strategies are proposed but others can be developed by using Theorem 1.

#### A. Using a threshold-like variable

First, we show that the event-triggering strategy proposed in [20] for sampled-data systems is not directly applicable to distributed NCS. Hence, we redesign this technique as in [17] by introducing an auxiliary variable. It has to be noted that the method in Section V.B in [17] cannot be applied 'off-the-shelf' here as we need to adapt the strategy to the protocol. We suppose that the controller (2) has been designed to make the closed-loop system (1) input-to-stable w.r.t. networked-induced error, which is equivalent to the following assumption (see Theorem 1 in [19]).

**Assumption 1.** *There exists a smooth Lyapunov function $V : \mathbb{R}^{n_x} \to \mathbb{R}$, $\underline{\alpha}_V, \overline{\alpha}_V, \alpha, \gamma \in \mathcal{K}_\infty$ such that for all $x \in \mathbb{R}^{n_x}$:*

$$\underline{\alpha}_V(|x|) \leq V(x) \leq \overline{\alpha}_V(|x|), \quad (8)$$

*and for all $(x, e) \in \mathbb{R}^{n_x+n_e}$:*

$$\langle \nabla V(x), f_x(x, e) \rangle \leq -\alpha(V(x)) + \gamma(|e|). \quad (9)$$

We suppose that the protocol is uniformly globally asymptotically stable (UGAS) [15], i.e., that the following holds.

**Assumption 2.** *There exist $W : \mathbb{R}_{\geq 0} \times \mathbb{R}^{n_e} \to \mathbb{R}_{\geq 0}$, $\underline{\alpha}_W, \overline{\alpha}_W \in \mathcal{K}_\infty$ and $\rho \in [0, 1)$ such that for all $(\kappa, e) \in \mathbb{Z}_{\geq 0} \times \mathbb{R}^{n_e}$ the following is satisfied:*

$$\underline{\alpha}_W(|e|) \leq W(\kappa, e) \leq \overline{\alpha}_W(|e|) \quad (10)$$
$$W(\kappa + 1, h_e(\kappa, e)) \leq \rho W(\kappa, e). \quad (11)$$

The round-robin and try-once-discard protocols have been shown to satisfy this property in [14], as well as other protocols (see [15] for instance). We note that when $\rho = 0$ we recover the situation in [17] where all nodes transmit at each transmission instant (i.e. $h_e(\kappa, e) = 0$ in (6)) (notice that (10) implies (11) in that case).

In view of (9), for any $\sigma \in \mathcal{K}_\infty$ with $\sigma(s) < s$ for $s > 0$, we have that $\gamma(|e|) \leq \sigma \circ \alpha(V(x))$ implies:

$$\langle \nabla V(x), f_x(x, e) \rangle \leq -(\mathbb{I} - \sigma) \circ \alpha(V(x)). \quad (12)$$

Instead of comparing $|e|$ and $|x|$ to derive the triggering rule as in [20], [17], we use the Lyapunov function $W(\kappa, e)$ which is characteristic of the protocol. According to (10), we have that $|e| \leq \underline{\alpha}_W^{-1}(W(\kappa, e))$. Therefore we can conclude that $\tilde{\gamma}(W(\kappa, e)) \leq V(x)$ (where $\tilde{\gamma}(s) = \alpha^{-1} \circ \sigma^{-1} \circ \gamma \circ \underline{\alpha}_W^{-1}(s)$, for $s \geq 0$) implies $\gamma(|e|) \leq \sigma \circ \alpha(V(x))$, that in

return ensures (12). Following the main idea of [20], a first attempt to define the triggering rule is:

$$\tilde{\gamma}(W(\kappa, e)) \geq V(x). \quad (13)$$

In this way, the flow and the jump sets of the corresponding system (6) are:

$$\begin{aligned} C &= \left\{(x, e, \kappa) : \tilde{\gamma}(W(\kappa, e)) \leq V(x)\right\} \\ D &= \left\{(x, e, \kappa) : \tilde{\gamma}(W(\kappa, e)) \geq V(x)\right\}. \end{aligned} \quad (14)$$

The problem with this policy is that we have no guarantee that $(x, e, \kappa)$ enters into $C$ after a jump. Indeed, while in [20] after each jump $e$ is reset to 0, here typically only a subvector $e_i$ is reset to zero after each transmission (see Section II). This may not be enough for $\tilde{\gamma}(W(\kappa, e))$ to become less than $V(x)$. As a consequence, the triggering rule (13) may generate several transmissions in a row before entering into $C$ that is unrealistic and contradicts item (iv) of Theorem 1. To overcome this drawback, we introduce a variable $\eta \in \mathbb{R}_{\geq 0}$ with the following dynamics:

$$\dot{\eta} = -\delta(\eta), \qquad \eta^+ = \tilde{\gamma}(W(\kappa, e)), \quad (15)$$

where $\delta$ is any locally Lipschitz class-$\mathcal{K}_\infty$ function. The system is now modeled as:

$$\left.\begin{aligned} \dot{x} &= f_x(x, e) \\ \dot{e} &= f_e(x, e) \\ \dot{\kappa} &= 0 \\ \dot{\eta} &= -\delta(\eta) \end{aligned}\right\} q \in C, \quad \left.\begin{aligned} x^+ &= x \\ e^+ &= h_e(\kappa, e) \\ \kappa^+ &= \kappa + 1 \\ \eta^+ &= \tilde{\gamma}(W(\kappa, e)) \end{aligned}\right\} q \in D, \quad (16)$$

where $q = (x, e, \kappa, \eta)$ and the sets $C$ and $D$ are defined as

$$\begin{aligned} C &= \left\{q : \max\{V(x), \eta\} \geq \tilde{\gamma}(W(\kappa, e)) \text{ and } \eta \geq 0\right\} \\ D &= \left\{q : \max\{V(x), \eta\} \leq \tilde{\gamma}(W(\kappa, e)) \text{ and } \eta \geq 0\right\}. \end{aligned} \quad (17)$$

The variable $\eta$ can be regarded as a decreasing threshold on $\tilde{\gamma}(W)$ in view of (15), and thus we enter into $C$ after a jump. Indeed, we have that $\eta^+ = \tilde{\gamma}(W(\kappa, e)) \geq \tilde{\gamma}(\rho W(\kappa, e)) \geq \tilde{\gamma}(W(\kappa^+, e^+))$ according to (11), (15) and since $\rho < 1$ and $\tilde{\gamma}$ is strictly increasing, therefore $q^+ \in C$. We are now able to apply Theorem 1 to guarantee stability properties for system (16) and the existence of dwell times.

**Theorem 2.** *Consider system (16) and suppose the following conditions hold.*

(i) *Assumptions 1-2 are satisfied.*
(ii) *Function $\tilde{\gamma}(W)$ is locally Lipschitz in $e$.*
(iii) *For any compact set $S \subset \mathbb{R}^{n_x + n_e}$, there exist $L_1, L_2 \in \mathbb{R}_{\geq 0}$ such that for any $(x, e) \in S$, $\kappa \in \mathbb{Z}_{\geq 0}$,*

$$\begin{aligned} |\langle \nabla V(x), f_x(x, e)\rangle| &\leq L_1(V(x) + \tilde{\gamma}(W(\kappa, e))) \\ |\tilde{\gamma}(W)^\circ(e; f_e(x, e))| &\leq L_2(V(x) + \tilde{\gamma}(W(\kappa, e))). \end{aligned}$$

(iv) *There exists $\epsilon \in [0, 1)$ such that $\lim_{s \to 0} \frac{\tilde{\gamma}(\rho s)}{\tilde{\gamma}(s)} = \epsilon$.*

*Then $\mathcal{A} = \{q : (x, e, \eta) = 0\}$ is S-GAS and solutions to (16) have a semiglobal dwell time on $\mathbb{R}^{n_q} \setminus \mathcal{A}$.*

Items (iii) and (iv) need to be added to guarantee the existence of dwell times compared to Theorem 3 in [17]. It is easy to check that condition (iv) of Theorem 2 holds when $\tilde{\gamma}$ is polynomial or homogeneous of degree $k \in \mathbb{Z}_{>0}$.

### B. Using a clock-like variable

In [6], NCS with time-triggered execution are modeled as a hybrid system similar to (6) by introducing a clock variable $\tau$ that would correspond to $\eta = \tau$ in (6). The flow and the jump sets are defined as $\tau$ being bigger or not than a given fixed bound $T$ known as MATI. This constant $T$ corresponds to the time it takes for the solution of the ordinary differential equation $\dot{\zeta} = -2L\zeta - \gamma(\zeta^2 + 1)$ to decrease from $\rho^{-1}$ to $\rho$, where $L$ and $\gamma$ are some constants (see (5) in [6]) and $\rho \in (0, 1)$ is given by Assumption 2, which is assumed to hold. In this subsection, we reduce the conservativeness of the strategy in [6] by making the ordinary differential equation that defines $\zeta$ state-dependent. This allows us to consider a larger class of systems and to potentially enlarge the inter-execution intervals compared to [6]. We suppose that Assumption 2 is satisfied with $W$ locally Lipschitz in $e$ and that the following holds.

**Assumption 3.** *There exist a locally Lipschitz function $V : \mathbb{R}^{n_x} \to \mathbb{R}$ and continuous functions $H : \mathbb{R}^{n_x} \to \mathbb{R}_{\geq 0}$, $L, G : \mathbb{R}^{n_x + n_e} \to \mathbb{R}_{\geq 0}$, $\underline{\alpha}_V, \overline{\alpha}_V \in \mathcal{K}_\infty$, $\varrho : \mathbb{R} \to \mathbb{R}$ continuous, positive definite such that the following conditions holds.*

(i) *For all $(x, e, \kappa) \in \mathbb{R}^{n_x + n_e} \times \mathbb{Z}_{\geq 0}$,*

$$W^\circ(e; f_e(x, e)) \leq L(x, e)W(\kappa, e) + H(x), \quad (18)$$

*where $W$ comes from Assumption 2.*

(ii) *For all $x \in \mathbb{R}^{n_x}$:*

$$\underline{\alpha}_V(|x|) \leq V(x) \leq \overline{\alpha}_V(|x|). \quad (19)$$

(iii) *For all $(x, e, \kappa) \in \mathbb{R}^{n_x + n_e} \times \mathbb{Z}_{\geq 0}$:*

$$\begin{aligned} V^\circ(x; f_x(x, e)) &\leq -\varrho(|x|) - \varrho(|e|) - H^2(x) \\ &\quad + G(x, e)W^2(\kappa, e). \end{aligned} \quad (20)$$

In [6], $L$ and $G$ are supposed to be constant that implies that system $\dot{x} = f_x(x, e)$ is $\mathcal{L}_2$-gain stable from $W$ to $H$. Making $L$ and $G$ state-dependent allows us to enlarge the studied class of systems and to eventually obtain less conservative upper bounds in (18) and (20) that will help to enlarge the inter-event intervals. Model (6) becomes here:

$$\left.\begin{aligned} \dot{x} &= f_x(x, e) \\ \dot{e} &= f_e(x, e) \\ \dot{\kappa} &= 0 \\ \dot{\eta} &= -2\eta L(x, e) - \eta^2 - G(x, e) \\ x^+ &= x \\ e^+ &= h_e(\kappa, e) \\ \kappa^+ &= \kappa + 1 \\ \eta^+ &= a \end{aligned}\right\} q \in C \quad q \in D, \quad (21)$$

where $q = (x, e, \kappa, \eta)$, $a$ is any constant in $(\rho, \infty)$ and $\eta \in \mathbb{R}$ plays the role of $\zeta$ mentioned above and is called a clock-like variable (see [6]). The sets $C$ and $D$ are:

$$C = \{q : \eta \in [a\rho^2, a]\}, \quad D = \{q : \eta = a\rho^2\}, \quad (22)$$

Note that, instead of setting $\eta^+$ to $\rho^{-1}$ at jumps as in [6], we consider any $a \in (\rho, \infty)$ that may help generating larger inter-execution intervals.

**Remark 1.** *We focus in this paper on NCS for which only one node among the $l$ nodes communicates at each transmission instant so that $\rho > 0$. When $\rho = 0$, as in [17],* we can redefine the sets in (22) as follows: $C = \{q : \eta \in [b, c]\}$, $D = \{q : \eta = b\}$ where $0 < b < c$ some constants.

The following theorem ensures the stability of system (21) and the existence of dwell times.

**Theorem 3.** *Consider system (21) and suppose Assumption 2-3 hold with $W$ locally Lipschitz in $e$. Then the set $\mathcal{A} = \{q : (x, e) = 0\}$ is S-GAS and solutions to (21) have a semiglobal dwell time on $\mathbb{R}^{n_q} \backslash \mathcal{A}$.*

## V. Self-triggered control

As mentioned before, the proposed event-triggering schemes in Section IV require the coordinator to continuously evaluate the triggering condition. This may induce a significant cost in terms of communication, computation time and power. To overcome these drawbacks, a possible solution lies in the self-triggered strategy. Self-triggered control considers the mathematical model of the control system and the last measurement of the plant states and/or the last control input in order to derive the next transmission instant. It represents a model-based emulation of event-triggered control in the sense that it identifies the time instants at which the jump condition is satisfied. In this section, we suppose that an event-triggering strategy has been designed such that the following conditions hold.

**Assumption 4.** *The conditions of Theorem 1 hold and item (iii) is satisfied on $C \cup D$ for system (6) with flow and jump sets of the form $C = \{q : \Gamma(q) \leq 0\}$ and $D = \{q : \Gamma(q) \geq 0\}$ where $\Gamma : \mathbb{R}^{n_q} \to \mathbb{R}_{\geq 0}$.*

Formally speaking, the event-triggered control strategy in Assumption 4 requires data transmission at the following time instant, for $j \in \mathbb{Z}_{>0}$:

$$\tau(j) = \inf\{t > t_{j-1}, (t, j-1) \in \text{dom}\,\phi : \Gamma(\phi(t, j-1)) = 0\},$$

where $\phi$ is a solution to (6). In order to guarantee stability, data transmission needs to occur no later than $\tau(j)$. When only the previous measurement of the state is available, the computation of $\tau(j)$ in an exact way is in most cases not possible. The self-triggered strategy computes a lower bound for $\tau(j)$ at which the next jump will occur. In [2], this lower bound is taken to be $\hat{\tau}(j) = \lambda(j)t_*$, where $\lambda(j)$ is strictly positive and satisfies for each $(t_j, j) \in \text{dom}\,\phi$:

$$\sum_{i=0}^{n-1} \varsigma_i (\mathcal{L}^i_{f_q}\Gamma)(\phi(t_j, j))\lambda^i(j) = 0, \tag{23}$$

where we have denoted the $i$th Lie derivative of $\Gamma$ along $f$ as $\mathcal{L}^i_f \Gamma = \mathcal{L}_f(\mathcal{L}^{i-1}_f \Gamma)$, $(\mathcal{L}_f \Gamma)(x) = \frac{\partial \Gamma}{\partial x}f(x)$ and $\mathcal{L}^0_f \Gamma = f$. The parameters $t_* \in \mathbb{R}_{>0}$ and $\varsigma_i \in \mathbb{R}$ are coefficients computed from $f_q$ and $\Gamma$. We assume that $\Gamma$ and $f_q$ are smooth functions in $(x, e, \eta)$. By abuse of notation we consider the $i$th Lie derivative $\mathcal{L}^i_{f_q}\Gamma$ even though $f_q$ and $\Gamma$ may not be differentiable in $\kappa$ (this is justified since $\dot{\kappa} = 0$ on flows). Equation (23) corresponds to the bound provided in Theorem V.4 in [2]. Guidelines for the design parameter $t_*$ are provided in Section VIII-A in [2], and the set of parameters $\varsigma_i$ have to be chosen to satisfy inequality (V.12) in [2]. We refer to [2] for a more detailed description of the roles played by these coefficients and how to compute them from the expressions of $f_q$ and $\Gamma$. The parameter $n > 1$ in (23) represents a design choice that trades the accuracy of the bound for the computational complexity. In other words, high values of $n$ imply times $\hat{\tau}(j)$ closer to $\tau(j)$ (and therefore less transmissions), but at the cost of solving a more complex algebraic equation. Since we may now transmit before $\Gamma(q) = 0$, we suppose item (iii) of Theorem 1 holds on $C \cup D$ (and not only $D$) in order to guarantee that the considered Lyapunov function still does not increase at jumps. This additional condition typically comes for free as it is the case in Section IV for instance. The problem is modeled as follows:

$$\left.\begin{array}{rcl}\dot{q} &=& f_q(q) \\ \dot{\tau}_1 &=& 1 \\ \dot{\tau}_2 &=& 0\end{array}\right\} \tilde{q} \in \widetilde{C} \\ \left.\begin{array}{rcl}q^+ &=& h_q(q) \\ \tau_1^+ &=& 0 \\ \tau_2^+ &=& \max\{\lambda(\kappa)t_*, \varepsilon\}\end{array}\right\} \tilde{q} \in \widetilde{D}, \tag{24}$$

where $\tau_1 \in \mathbb{R}_{\geq 0}$ is a clock variable and $\tau_2 \in \mathbb{R}_{\geq 0}$ is used to define the next transmission instant, $\tilde{q} = (q, \tau_1, \tau_2)$, and the sets $\widetilde{C}$ and $\widetilde{D}$ are:

$$\widetilde{C} = \left\{\tilde{q} : \tau_1 \in [0, \tau_2]\right\}, \quad \widetilde{D} = \left\{\tilde{q} : \tau_1 \geq \tau_2\right\}. \tag{25}$$

It is shown in [2] that the formula in (23) can be used to design $\hat{\tau}(j)$ that is very close to $\tau(j)$. Nevertheless, to prevent from the situation where the self-triggering technique of [2] generates conservative times because of an inadequate parameters choice in (23), we introduce an arbitrary small constant $\varepsilon \in \mathbb{R}_{>0}$ in (24) to guarantee the existence of a minimal interval of time between two transmissions. This is justified since Assumption 4 implies the existence of such a constant time (semiglobally), in view of item (iv) of Theorem 1. The following theorem shows that the properties ensured by the considered event-triggering strategy are maintained under the proposed self-triggering rules.

**Theorem 4.** *Let $f_q$ and $\Gamma$ be smooth functions in $(x, e, \eta)$. Under Assumption 4, the set $\widetilde{\mathcal{A}} = \{\tilde{q} : (x, e, \upsilon(\eta)) = 0\}$ is S-GAS for system (24) and solutions to (24) have a semiglobal dwell time on $\mathbb{R}^{n_{\tilde{q}}} \backslash \widetilde{\mathcal{A}}$.*

Self-triggering rules can be derived for the event-triggered strategies developed in Section IV as follows. For Section IV-B, we take $\Gamma(q) = \eta - a\rho^2$. The function $\Gamma$ is smooth and we assume that $f_x, f_e, L, G$ are smooth functions in $(x, e, \eta)$. Embedding system (21) into (24) allows us to obtain a self-triggered control technique. The conclusions of Theorem 4 apply as all the required conditions hold. For the event-triggered control of Section IV-A, we cannot define $\Gamma(q) = \tilde{\gamma}(W(\kappa, e)) - \max\{V(x), \eta\}$ as $\Gamma$ will not be smooth. Therefore, we define $\Gamma_1(q) = \tilde{\gamma}(W(\kappa, e)) - \eta$ and $\Gamma_2(q) = \tilde{\gamma}(W(\kappa, e)) - V(x)$. By assuming that $\tilde{\gamma}(W)$ and $V$ are smooth in $e$ and $x$ respectively, we see that $\Gamma_1$ and $\Gamma_2$ are smooth in $(x, e, \eta)$. We modify the jump equation for $\tau_2$ in (24) as follows: $\tau_2^+ = \max\{\lambda_1(\kappa)t_*, \lambda_2(\kappa)t_*\}$ where $\lambda_i$ satisfies (23) with $\Gamma = \Gamma_i$, $i \in \{1, 2\}$. Under this setup, the conditions of Theorem 4 are satisfied provided that $f_x, f_e, \delta$ are smooth.

| n | $\varsigma_0$ | $\varsigma_1$ | $\varsigma_2$ | $\varsigma_3$ |
|---|---|---|---|---|
| $\Gamma_1$: 4 | $-8.06 \cdot 10^3$ | $-226.07$ | $481.76$ | $-258.47$ |
| $\Gamma_2$: 3 | $-1.46 \cdot 10^3$ | $-1.21 \cdot 10^3$ | $4.94 \cdot 10^3$ | - |

TABLE I

PARAMETERS FOR THE SELF-TRIGGERED TECHNIQUE.

## VI. ILLUSTRATIVE EXAMPLE

To illustrate the proposed strategy, we consider the control of a jet engine compressor. We borrow the model from [10]:

$$\dot{x}_1 = -x_2 - \frac{3}{2}x_1^2 - \frac{1}{2}x_1^3, \qquad \dot{x}_2 = u, \qquad (26)$$

where $x_1$ represents the mass flow, $x_2$ is the pressure rise and $u$ is the throttle mass flow. In this model the origin has been translated to the desired equilibrium point, hence the objective is to steer $(x_1, x_2)$ to zero. The control law $u = 4x_1 - 4x_2 - \frac{9}{2}x_1^2 - \frac{3}{2}x_1^3$ is designed to stabilize the system. This controller is connected to the two sensors measuring $x_1$ and $x_2$ through a network under the TOD protocol. For simplicity, we consider that the controller is connected to the actuator. Nonetheless, as pointed out previously in this paper, the developed framework accounts for the more realistic case of the controller and the actuator not being collocated. We implement the event-triggering strategy proposed in Section IV-A. Assumption 1 is satisfied with $V(x) = \frac{1}{2}x_1^2 + \frac{1}{2}(x_2 - 3x_1)^2$, $\alpha(s) = -0.066s$, $\gamma(s) = 4.37 \cdot 10^4 s^2 + 9.10 \cdot 10^6 s^4$ for $s \geq 0$. Assumption 2 is satisfied for $W(e) = |e|$ for the TOD protocol, where $e_1$ and $e_2$ represent the network-induced errors for $x_1$ and $x_2$, and with $\rho = \sqrt{\frac{l-1}{l}} = \frac{1}{\sqrt{2}}$ as the number of nodes is $l = 2$ (see Proposition 5 in [14]). We select $\sigma(s) = 0.9s$ in (12) and thus we obtain $\tilde{\gamma}(W(\kappa, e)) := 7.34 \cdot 10^5 |e|^2 + 1.52 \cdot 10^8 |e|^4$. The Yalmip software [12] was used to compute $\alpha$, $\gamma$ and $\tilde{\gamma}(W(\kappa, e))$. We now construct a linear differential equation for the auxiliary variable $\eta$: $\dot{\eta} = -0.01\eta$, with initial condition $\eta(0, 0) = 5000$. It is expected that larger values of $\eta(0, 0)$ would enlarge the transmission times, at the cost of a degradation in performance. It can be verified that items (ii)-(iv) of Theorem 2 hold (we use the fact that $\tilde{\gamma}$ is a polynomial function to show item (iv) of Theorem 2). We derive as well a self-triggered emulation of the event-triggered technique of Section IV-A. As detailed at the end of Section V, we consider $\Gamma_1(q) := \tilde{\gamma}(W(\kappa, e)) - \eta$ and $\Gamma_2(q) := \tilde{\gamma}(W(\kappa, e)) - V(x)$. The design parameters for the self-triggering condition are reported in Table I. In both cases we considered $t_* = 10^{-3}$.

We compare the event-triggered strategy herein proposed with a periodic implementation. In order to compute a period, we apply the technique in [6]. For an operating ball of radius 1, the obtained period is $T = 0.010$. For the comparison, we consider 200 different initial conditions randomly distributed in a ball of radius 1. Table II shows the average inter-transmission time under the three different strategies. Both the event-triggered and the self-triggered strategy outperform the periodic approach. The gap between the event-triggered and the self-triggered inter-transmission times is due to the conservativeness of the technique in [2]. To further illustrate the proposed approach, we depict as well the evolution of the transmission times and the network-induced error under

| Periodic [6] | Event-triggered | Self-triggered |
|---|---|---|
| 0.010 | 0.061 | 0.046 |

TABLE II

AVERAGE INTER-TRANSMISSION TIME FOR 200 INITIAL CONDITIONS.

the self-triggered strategy for a particular initial condition $(x(0, 0) = (0.95, -0.14), e(0, 0) = (0, 0))$. The plot shows how the network grants access to the node with the largest error. Likewise, the transmission times vary according to the current state of the plant. This fact suggests that the rigid periodic paradigm overloads unnecessarily the network, and the flexibility and adaptability of event-triggered control (and consequently self-triggered control) is able to relax this requirement and reduce network usage.

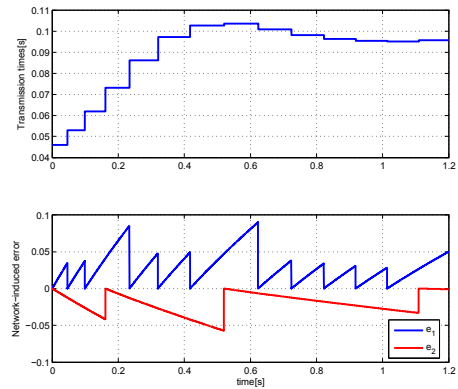

Fig. 1. Evolution of the transmission times and the networked-induced error under the self-triggered strategy.

## VII. CONCLUSION

We have used the prescriptive framework in [17] originally developed for sampled-data systems to synthesize novel event-triggering rules for distributed NCS. We have then shown how self-triggering conditions can then been derived (under some conditions) by applying the techniques in [2]. This work represents the first step towards a more fundamental question in distributed networked control: given a set of physically distributed sensors and actuators, how should communication between the different nodes be scheduled? Ad-hoc solutions to this problem include the TOD protocol, where the node with the largest network-induced error information is granted access to the communication channel. The presented framework can be further explored to design communication protocols that decide the order at which each node needs to send information.

## VIII. APPENDIX

**Proof of Theorem 1.** Let $B = \{q \in C \cup D : |q|_\mathcal{A} < \Delta\}$ where $\Delta \in \mathbb{R}_{>0}$ and $\phi$ be a solution to (6) with $\phi(0, 0) \in B$. Define the set $\Omega = \{q \in C \cup D : R(q) \leq \mu\}$ where $\mu \in \mathbb{R}_{>0}$ is such that $B \subseteq \Omega$ (take for instance $\mu = \overline{\alpha}_R(\Delta)$ in view of item (i) of Theorem 1). The set $C \cup D$ is forward invariant for system (6) since $G(D) \subset (C \cup D)$. Hence, in view of items (i)-(iii) of Theorem 1, $\phi(t, j) \in \Omega$ for any $(t, j) \in \text{dom} \, \phi$. From item (ii) of Theorem 1 and by using standard comparison principles, there exists $\beta \in \mathcal{KL}$ that satisfies, for all $(s, t_1, t_2) \in \mathbb{R}_{\geq 0} \times \mathbb{R}_{\geq 0} \times \mathbb{R}_{\geq 0}$:
$$\beta(s, t_1 + t_2) = \beta(\beta(s, t_1), t_2), \quad (27)$$
and such that, for all $(t_j, j) \preceq (t, j) \in \operatorname{dom} \phi$,
$$R(\phi(t, j)) \leq \beta(R(\phi(t_j, j)), t - t_j), \quad (28)$$
where $(t_j, j) \preceq (t, j)$ means that $t_j \leq t$. From item (iii) of Theorem 1, it follows that:
$$R(\phi(t_{j+1}, j+1)) \leq R(\phi(t_{j+1}, j)) \quad (29)$$
for all $j$ such that $(t, j) \in \operatorname{dom} \phi$ for some $t \in \mathbb{R}_{\geq 0}$. Combining (27)-(29), we obtain:
$$R(\phi(t, j)) \leq \beta(R(\phi(0, 0)), t) \quad \forall (t, j) \in \operatorname{dom} \phi. \quad (30)$$

Now let $(t, j) \in \operatorname{dom} \phi$, if $\phi(t, j) \in \mathcal{A}$ then $R(\phi(t, j)) = 0$ according to item (i) of Theorem 1. If $\phi(t, j) \notin \mathcal{A}$ then that means that $\phi(t', j') \notin \mathcal{A}$ for all $(t', j') \preceq (t, j) \in \operatorname{dom} \phi$ (i.e. $t' \leq t$ and $j' \leq j$) since $\mathcal{A}$ is forward invariant for system (6) (according to items (i)-(iii) of Theorem 1 and since $G(D) \subset (C \cup D)$). As a consequence, we have that $t' \geq \tau(\mu) j'$ for all $(t', j') \preceq (t, j) \in \operatorname{dom} \phi$, where $\tau(\Delta) \in \mathbb{R}_{>0}$ is a minimal interval of times between two jumps on $\Omega \backslash \mathcal{A}$ whose existence is ensured by item (iv) of Theorem 1. It follows that:
$$R(\phi(t, j)) \leq \beta(R(\phi(0, 0)), \tfrac{1}{2} t + \tfrac{1}{2} \tau(\Delta) j). \quad (31)$$
By using item (i) of Theorem 1, we deduce that, for all $(t, j) \in \operatorname{dom} \phi$:
$$|\phi(t, j)|_{\mathcal{A}} \leq \underline{\alpha}_R^{-1}\Big(\beta\big(\overline{\alpha}_R(|\phi(0, 0)|_{\mathcal{A}}), \tfrac{1}{2} t + \tfrac{1}{2} \tau(\Delta) j\big)\Big), \quad (32)$$
denoting $\beta_\Delta : (s, t, j) \mapsto \underline{\alpha}_R^{-1}\Big(\beta\big(\overline{\alpha}_R(s), \tfrac{1}{2} t + \tfrac{1}{2} \tau(\Delta) j\big)\Big) \in \mathcal{KLL}$ (since $\beta \in \mathcal{KL}$, $\underline{\alpha}_R, \overline{\alpha}_R \in \mathcal{K}_\infty$), for all $(t, j) \in \operatorname{dom} \phi$:
$$|\phi(t, j)|_{\mathcal{A}} \leq \beta_\Delta(|\phi(0, 0)|_{\mathcal{A}}, t, j). \quad (33)$$
Hence, the set $\mathcal{A}$ is S-GAS according to Definition 1 in [17]. $\square$

**Proof of Lemma 1.** Let $\Delta \in \mathbb{R}_{>0}$ and define $B = \{q \in C \cup D : |q|_{\mathcal{A}} < \Delta\}$ and $\Omega = \{q \in C \cup D : R(q) \leq \mu\}$ where $\mu \in \mathbb{R}_{>0}$ is such that $B \subseteq \Omega$ (take for instance $\mu = \overline{\alpha}_R(\Delta)$ in view of item (i) of Theorem 1). According to items (i)-(iii) of Theorem 1, $\phi(t, j) \in \Omega$ for any $(t, j) \in \operatorname{dom} q$. We notice that $\Omega = \mathcal{A} \cup \Theta(\mu)$. Denote $t_1 \in (0, \infty)$ the first jump instant (if no jump ever occurs i.e. $t_1 = \infty$, (2) in [17] is obviously satisfied). We have that $\phi(t_1, 0) \in D$ and $R(\phi(t_1, 1)) \leq \mu$. When $R(\phi(t_1, 1)) = 0$ then $\phi(t, j) \in \mathcal{A}$ for any $(t_1, 1) \preceq (t, j) \in \operatorname{dom} q$ (i.e. $t_1 \leq t$ and $1 \leq j$) as $\mathcal{A}$ is forward invariant for system (6) and (2) in [17] is ensured. When $\phi(t_1, 1) \in \Theta(\mu)$, then according to item (i) of Lemma 1, $\psi(h_q(\phi(t_1, 0))) = \psi(\phi(t_1, 1)) \leq a < b$ and therefore a jump cannot occur immediately in view of item (iii-b) of Lemma 1. Let denote $t_2 > t_1$ with $(t_2, 2) \in \operatorname{dom} \phi$ the next jump instant and suppose $t_2 \neq \infty$, otherwise the desired result holds. By the continuity of $\psi$ and the solution $\phi$ to system (6) on flows, there exists $t^* \in (t_1, t_2]$ such that $\psi(\phi(t^*, 1)) = b$ for any $t \in [t_1, t^*]$ with $(t, 1) \in \operatorname{dom} \phi$ as $a < b$. According to item (ii) of Lemma 1, we deduce that $\phi(t, 1) \in C \backslash D$ for any $t \in [t_1, t_2)$. In view of item (iii) of Lemma 1, invoking standard comparison principles, we deduce that $\psi(\phi(t, 1)) \leq \theta(t)$ for any $t \in [t_1, t^*]$ where $\theta(t)$ is the solution of $\dot{\theta} = \lambda(\theta)$ satisfying $\theta(t_1) = a \geq \psi(\phi(t_1, 1))$. The next jump cannot occur before the time $\tau(\Delta) \in (0, \infty]$ it takes for $\theta$ to evolve from $a$ to $b$ (which is independent of $t_1$) has elapsed. By induction, we deduce that the inter-jump interval on $B \backslash \mathcal{A}$ is lower bounded by $\tau(\Delta)$. Hence solutions to (6) have a semiglobal dwell time on $\mathbb{R}^{n_q} \backslash \mathcal{A}$ according to Definition 2 in [17]. $\square$

**Proof of Theorem 2.** We verify the conditions of Theorem 1 and apply it to obtain the desired results. We note that $G(D) = \{q : \kappa \in \mathbb{Z}_{>0} \text{ and } \eta \geq 0\} \subset C \cup D$. We consider the following candidate Lyapunov function:
$$R(q) = \max\{V(x), \tilde{\gamma}(W(\kappa, e)), \eta\}. \quad (34)$$
It can be shown that item (i) of Theorem 1 holds with $\upsilon(\eta) = \eta$ by using (8), (10) and Remark 2.3 in [11] and noting that $\eta \geq 0$ on $C \cup D$. On $C$, we have that $\tilde{\gamma}(W(\kappa, e)) \leq \max\{V(x), \eta\} = R(q)$, as a consequence, when $R(q) = V(x)$, we get in view of (12): $R^\circ(q; f_q(q)) \leq -(\mathbb{I} - \sigma) \circ \alpha(V(x)) = -(\mathbb{I} - \sigma) \circ \alpha(R(q))$, and when $R(q) = \eta$, $R^\circ(q; f_q(q)) = -\delta(\eta) = -\delta(R(q))$. Thus, item (ii) of Theorem 1 is ensured with $\alpha_R(s) = \min\{(\mathbb{I} - \sigma) \circ \alpha(s), \delta(s)\}$ for $s \geq 0$. Let $q \in D$,
$$R(q^+) = \max\{V(x^+), \tilde{\gamma}(W(\kappa^+, e^+)), \eta^+\}, \quad (35)$$
according to Assumption 2,
$$\begin{aligned} R(q^+) &\leq \max\{V(x), \tilde{\gamma}(\rho W(\kappa, e)), \tilde{\gamma}(W(\kappa, e))\} \\ &\leq \max\{V(x), \tilde{\gamma}(W(\kappa, e))\} = R(q) \end{aligned} \quad (36)$$
since $\rho < 1$ and $\tilde{\gamma}$ is increasing: item (iii) of Theorem 1 is ensured. We now prove that item (iv) of Theorem 1 is satisfied using Lemma 1. Let $\psi : q \mapsto \frac{\tilde{\gamma}(W(\kappa, e))}{R(q)}$ that is defined on $\Theta(\mu)$ with $\mu \in \mathbb{R}_{>0}$ ($\Theta(\mu)$ comes from Lemma 1). Let $q \in D$ such that $h_q(q) \in \Theta(\mu)$, using Assumption 2 we have:
$$\psi(h_q(q)) = \frac{\tilde{\gamma}(W(\kappa+1, h_e(\kappa, e)))}{R(h_q(q))} \leq \frac{\tilde{\gamma}(\rho W(\kappa, e))}{R(h_q(q))} \quad (37)$$
since $R(h_q(q)) \geq \tilde{\gamma}(W(\kappa, e))$,
$$\psi(h_q(q)) \leq \frac{\tilde{\gamma}(\rho W(\kappa, e))}{\tilde{\gamma}(W(\kappa, e))}. \quad (38)$$
Define $\theta(\mu) = \max_{s \in (0, \tilde{\gamma}^{-1}(\mu)]} \frac{\tilde{\gamma}(\rho s)}{\tilde{\gamma}(s)}$. Function $s \mapsto \frac{\tilde{\gamma}(\rho s)}{\tilde{\gamma}(s)}$ is defined and continuous on $(0, \tilde{\gamma}^{-1}(\mu)]$ since $\tilde{\gamma}$ is continuous and only cancels at the origin. Moreover, $\frac{\tilde{\gamma}(\rho s)}{\tilde{\gamma}(s)} < 1$ for any $s \in (0, \tilde{\gamma}^{-1}(\mu)]$ since $\tilde{\gamma}$ is strictly increasing as a class-$\mathcal{K}_\infty$ function and $\rho < 1$. On the other hand, in view of item (iv) of Theorem 2, $\lim_{s \to 0} \frac{\tilde{\gamma}(\rho s)}{\tilde{\gamma}(s)} = \epsilon < 1$, therefore we can write that $\theta(\mu) < 1$. In that way, from (38), the following holds:
$$\psi(h_q(q)) \leq \theta(\mu), \quad (39)$$
and item (i) of Lemma 1 is satisfied with $a = \theta(\mu)$. Since after each jump, $\psi(h_q(q)) \leq \theta(\mu) < 1$ for $q \in \Theta(\mu)$, we deduce from the definition of $\psi$ that item (ii) of Lemma 1 holds with $b = 1$, since as long as $\tilde{\gamma}(W(\kappa, e)) < R(q)$, $q \in C \backslash D$. Let $q \in \Theta(\mu)$, we have that:

$$|R^\circ(q; f_q(q))| \leq \max\{\left|\frac{\partial V}{\partial x} f_x(x, e)\right|, \delta(\eta), |\tilde{\gamma}(W(\kappa, e))^\circ(q; f_q(q))|\} \quad (40)$$

in view of item (iii) of Theorem 2 and denoting $k$ the Lipschitz constant of $\delta$ on $\Theta(\mu)$ and $L = \max\{L_1, L_2\}$, we obtain:

$$|R^\circ(q; f_q(q))| \leq \max\{L(V(x) + \tilde{\gamma}(W(\kappa, e))), k|\eta|\} \quad (41)$$

applying the definition of $R$ in (34),

$$|R^\circ(q; f_q(q))| \leq \max\{L(R(q) + R(q)), kR(q)\} = \max\{2L, k\}R(q), \quad (42)$$

hence

$$|\tilde{\gamma}(W(\kappa, e))R^\circ(q; f_q(q))| \leq \max\{2L, k\}R^2(q). \quad (43)$$

On the other hand, from item (iii) of Theorem 2,

$$|\tilde{\gamma}(W(\kappa, e))^\circ(q; f_q(q))R(q)| \leq L_2(V(x) + \tilde{\gamma}(W(\kappa, e)))R(q) \leq 2L_2 R^2(q). \quad (44)$$

Consequently, from (40), (43) and (44):

$$\begin{aligned}\psi^\circ(q; f_q(q)) &= \frac{\tilde{\gamma}(W)^\circ(e; f_e(x,e))R(q) - R^\circ(q; f_q(q))\tilde{\gamma}(W(\kappa, e))}{R(q)^2} \\ &\leq \frac{2L_2 R^2(q) + \max\{2L, k\}R^2(q)}{R^2(q)} \\ &= 2L_2 + \max\{2L, k\} =: \lambda(\psi)\end{aligned} \quad (45)$$

where $\lambda$ is continuous and non-decreasing. We have proved that item (iii) of Lemma 1 applies: solutions to (16) have a semiglobal dwell time on $\mathbb{R}^{n_q} \backslash \mathcal{A}$. As a consequence, the set $\mathcal{A}$ is S-GAS according to Theorem 1. $\square$

**Proof of Theorem 3.** We show that the conditions of Theorem 1 holds. We note that $G(D) = \{q : \kappa \in \mathbb{Z}_{>0} \text{ and } \eta = a\} \subset \{q : \kappa \in \mathbb{Z}_{>0} \text{ and } \eta \in [a\rho^2, a]\} = C \cup D$. We consider the candidate Lyapunov function:

$$R(q) = V(x) + \eta W^2(\kappa, e), \quad (46)$$

that satisfies item (i) of Theorem 1 with $\overline{\alpha}_R(s) = \overline{\alpha}_V(s) + a\overline{\alpha}_W^2(s)$ and $\underline{\alpha}_R(s) = \min\{\underline{\alpha}_V(\frac{s}{2}), a\rho^2 \underline{\alpha}_W(\frac{s}{2})\}$ for $s \geq 0$ by identifying $\chi(\eta) = \eta$ and $\upsilon(\eta) = 0$ and using (10), (19), the fact that $\eta \in [a\rho^2, a]$ on $C \cup D$ and Remark 2.3 in [11]. For any $q \in C$, we have that:

$$\begin{aligned}R^\circ(q; f_q(q)) &= V^\circ(x; f_x(x, e)) + 2\eta W(\kappa, e)W^\circ(e; f_e(x, e)) \\ &\quad + (-2\eta L(x, e) - \eta^2 - G(x, e))W^2(\kappa, e),\end{aligned}$$

from items (i) and (iii) of Assumption 3,

$$\begin{aligned}R^\circ(q; f_q(q)) &\leq -\varrho(|x|) - \varrho(|e|) - H^2(x) + G(x, e)W^2(\kappa, e) \\ &\quad + (-2\eta L(x, e) - \eta^2 - G(x, e))W^2(\kappa, e) \\ &\quad + 2\eta W(\kappa, e)(L(x, e)W(\kappa, e) + H(x)) \\ &= -\varrho(|x|) - \varrho(|e|) - H^2(x) - \eta^2 W^2(\kappa, e) \\ &\quad + 2\eta W(\kappa, e)H(x),\end{aligned}$$

using the fact that $2\eta W(\kappa, e)H(x) \leq \eta^2 W^2(\kappa, e) + H^2(x)$, we get:

$$R^\circ(q; f_q(q)) \leq -\varrho(|x|) - \varrho(|e|) \leq -\rho(|(x, e)|), \quad (47)$$

consequently, item (ii) of Theorem 1 is satisfied with $\alpha_R(s) = \rho \circ \overline{\alpha}_R^{-1}(s)$ for $s \geq 0$. Let $q \in D$,

$$\begin{aligned}R(q^+) &= V(x^+) + \eta^+ W^2(\kappa^+, e^+) \\ &= V(x) + aW^2(\kappa^+, e^+).\end{aligned} \quad (48)$$

We obtain using (11):

$$R(q^+) \leq V(x) + a\rho^2 W^2(\kappa, e), \quad (49)$$

since $\eta = a\rho^2$ on $D$, we have:

$$\begin{aligned}R(q^+) &\leq V(x) + a\rho^2 W(\kappa, e) \\ &= V(x) + \eta W(\kappa, e) = R(q),\end{aligned} \quad (50)$$

and item (iii) of Theorem 1 is ensured. Finally, we apply Lemma 1 to prove the existence of dwell times. Let $\psi : q \mapsto \frac{a\rho^2}{\eta}$ that is defined on $\Theta(\mu)$ with $\mu \in \mathbb{R}_{>0}$ (see Lemma 1). For $q \in D$, we have that $\psi(h_q(q)) = \rho^2$ so item (i) of Lemma 1 holds with $a = \rho^2$. Moreover, it can be seen that item (ii) of Lemma 1 is satisfied with $b = 1$. Let $q \in \Theta(\mu)$,

$$\begin{aligned}\psi^\circ(q; f_q(q)) &= -a\rho^2 \frac{\dot{\eta}}{\eta^2} = \frac{a\rho^2}{\eta^2}(2\eta L(x, e) + \eta^2 + G(x, e)) \\ &= \frac{a\rho^2}{\eta} 2L(x, e) + a\rho^2 + \frac{a\rho^2}{\eta^2} G(x, e) \\ &= 2\psi L(x, e) + a\rho^2 + \frac{\psi^2}{a\rho^2} G(x, e)\end{aligned}$$

denoting $M = \max_{R(q) \leq \mu}\{2L(x, e), G(x, e)\}$ (that is well defined since $L$ and $G$ are continuous),

$$\psi^\circ(q; f_q(q)) \leq M\psi + a\rho^2 + M\frac{\psi}{a\rho^2} =: \lambda(\psi). \quad (51)$$

We see that $\lambda$ is non-decreasing and continuous on $\mathbb{R}_{>0}$. As a consequence, item (iii) of Lemma 1 holds that implies that item (iv) of Theorem 1 is guaranteed: solutions to (21) have a semiglobal dwell time on $\mathbb{R}^{n_q} \backslash \mathcal{A}$. We obtain the desired result by applying Theorem 1. $\square$

**Proof of Theorem 4.** Consider system (24) and the candidate Lyapunov function $\widetilde{R}(\widetilde{q}) = R(q)$ for $\widetilde{q} = (q, \tau_1, \tau_2) \in \mathbb{R}^{n_{\widetilde{q}}}$, where $R$ comes from Theorem 1. We have that $\widetilde{C} \cup \widetilde{D} = C \cup D \times \mathbb{R}_{\geq 0}^2$ and $|\widetilde{q}|_{\widetilde{\mathcal{A}}} = |q|_{\mathcal{A}}$, therefore $\underline{\alpha}_R(|\widetilde{q}|_{\widetilde{\mathcal{A}}}) \leq \widetilde{R}(\widetilde{q}) \leq \overline{\alpha}_R(|\widetilde{q}|_{\widetilde{\mathcal{A}}})$ for any $\widetilde{q} \in \widetilde{C} \cup \widetilde{D}$ in view item (i) of Theorem 1. According to [2], $\tau_2$ in (24) is such that for any $\tau_1 \leq \tau_2$, $\Gamma(q) \leq 0$. Hence, $\widetilde{q} \in \widetilde{C}$ implies $q \in C$. Consequently, we have that $\widetilde{R}^\circ(\widetilde{q}; \widetilde{f}_q(\widetilde{q})) \leq -\alpha_R(\widetilde{R}(\widetilde{q}))$ in view of item (ii) of Theorem 1. Item (iii) of Theorem 1 is assumed to hold on $C \cup D$, therefore we have that $\widetilde{R}(\widetilde{h}_q(\widetilde{q})) \leq \widetilde{R}(\widetilde{q})$. Solutions to (24) obviously have a semiglobal dwell time on $\mathbb{R}^{n_{\widetilde{q}}} \backslash \widetilde{\mathcal{A}}$ since a jump cannot occur on $\mathbb{R}^{n_{\widetilde{q}}} \backslash \widetilde{\mathcal{A}}$ after $\varepsilon$ seconds has elapsed which is a semiglobal dwell-time exhibited by the event-triggered strategy. By following the similar lines as in the proof of Theorem 1 and noting that $G(\widetilde{D}) = G(D) \times \{(\tau_1, \tau_2) : \tau_1 = 0\} \subset \widetilde{C} \cup \widetilde{D}$, (since $G(D) \subset (C \cup D)$ in view of Assumption 4), the set $\widetilde{\mathcal{A}}$ is S-GAS. $\square$

REFERENCES

[1] A. Anta and P. Tabuada. To sample or not to sample: self-triggered control for nonlinear systems. *IEEE Transactions on Automatic Control*, 55(9):2030–2042, 2010.
[2] A. Anta and P. Tabuada. Exploiting isochrony in self-triggered control. *Provisionally accepted for publication. arXiv 1009.5208*, 2011.

[3] K.E. Arzén. A simple event-based PID controller. In *14$^{th}$ IFAC World Congress, Beijing, China*, 1999.

[4] K.J. Aström and B.M. Bernhardsson. Comparison of Riemann and Lebesgue sampling for first order stochastic systems. In *CDC (IEEE Conference on Decision and Control), Las Vegas, U.S.A.*, 2002.

[5] C. Cai, A.R. Teel, and R. Goebel. Smooth Lyapunov functions for hybrid systems part II: (pre)asymptotically stable compact sets. *IEEE Transactions on Automatic Control*, 53(3):734–748, 2008.

[6] D. Carnevale, A.R. Teel, and D. Nešić. A Lyapunov proof of an improved maximum allowable transfer interval for networked control systems. *IEEE Transactions on Automatic Control*, 52(5), 2007.

[7] M.C.F. Donkers and W.P.M.H. Heemels. Output-based event-triggered control with guaranteed $\mathcal{L}_\infty$-gain and improved event-triggering. In *IEEE Conference on Decision and Control, Atlanta*, 2010.

[8] R. Goebel and A.R. Teel. Solution to hybrid inclusions via set and graphical convergence with stability theory applications. *Automatica*, 42:573–587, 2006.

[9] W.P.M.H. Heemels, J.H. Sandee, and P.P.J. van den Bosch. Analysis of event-driven controllers for linear systems. *International Journal of Control*, 81(4):571–590, 2009.

[10] M. Krstić and P. Kokotović. Lean backstepping design for a jet engine compressor model. *4th IEEE Conf. on Control Applications*, 1995.

[11] D.S. Laila and D. Nešić. Lyapunov based small-gain theorem for parameterized discrete-time interconnected ISS systems. In *CDC (IEEE Conference on Decision and Control), Las Vegas, U.S.A.*, pages 2292–2297, 2002.

[12] J. Löfberg. Pre- and post-processing sum-of-squares programs in practice. *IEEE Transactions on Automatic Control*, 54(5), 2009.

[13] M. Mazo and P. Tabuada. Decentralized event-triggered control over wireless sensor/actuator networks. *Accepted for publication in IEEE Transactions on Automatic Control. arXiv:1004.0477*.

[14] D. Nešić and A.R. Teel. Input-output stability properties of networked control systems. *IEEE Transactions on Automatic Control*, 2004.

[15] D. Nešić and A.R. Teel. Input-to-state stability of networked control systems. *Automatica*, 40:2121–2128, 2004.

[16] D. Nešić, L. Zaccarian, and A.R. Teel. Stability properties of reset systems. *Automatica*, 44:2019–2026, 2008.

[17] R. Postoyan, A. Anta, D. Nešić, and P. Tabuada. A unifying Lyapunov-based framework for the event-triggered control of nonlinear systems. In *IEEE Conference on Decision and Control and European Control Conference, Orlando, U.S.A.*, 2011.

[18] R. Postoyan, P. Tabuada, D. Nešić, and A. Anta. Event-triggered and self-triggered stabilization of distributed networked control systems. In *IEEE Conference on Decision and Control and European Control Conference, Orlando, U.S.A.*, 2011.

[19] E.D. Sontag and Y. Wang. On characterizations of the input-to-state stability property. *Systems & Control Letters*, 24(5):351–359, 1995.

[20] P. Tabuada. Event-triggered real-time scheduling of stabilizing control tasks. *IEEE Transactions on Automatic Control*, 52(9), 2007.

[21] G.C. Walsh, O. Beldiman, and L.G. Bushnell. Asymptotic behavior of nonlinear networked control systems. *IEEE Transactions on Automatic Control*, 46:1093–1097, 2001.

[22] X. Wang and M.D. Lemmon. Self-triggered feedback control systems with finite-gain $\mathcal{L}_2$ stability. *IEEE Transactions on Automatic Control*, 45:452–467, 2009.

[23] X. Wang and M.D. Lemmon. Asymptotic stability in distributed event-triggered networked control systems with delays. In *ACC (American Control Conference), Baltimore, U.S.A.*, 2010.